\documentclass[12pt]{amsart}
\usepackage{amssymb,amsmath}
\usepackage{amscd}
\usepackage{verbatim}
\makeatletter
\def\@cite#1#2{{\m@th\upshape\bfseries%
[{#1\if@tempswa{\m@th\upshape\mdseries, #2}\fi}]}}
\makeatother
\theoremstyle{plain}
\newtheorem{thm}{Theorem}[section]
\newtheorem{lem}[thm]{Lemma}
\newtheorem{cor}[thm]{Corollary}
\newtheorem{prop}[thm]{Proposition}
\theoremstyle{definition}
\newtheorem{rem}[thm]{Remark}

\newtheorem{defn}[thm]{Definition}
\newtheorem{note}[thm]{Note}

\newcommand{\Prf}{\noindent\textbf{Proof.\ }}
\newcommand{\bx}{\strut\hfill$\blacksquare$\medbreak}

\newcommand{\ca}{\mathrm{C}^*}

\newcommand{\ol}{\overline}

\newcommand{\td}{\widetilde}

\newcommand{\bbT}{{\mathbb{T}}}
\newcommand{\bbZ}{{\mathbb{Z}}}
 \newcommand{\A}{{\mathcal{A}}}
 \newcommand{\B}{{\mathcal{B}}}
 \newcommand{\C}{{\mathcal{C}}}
 
 \newcommand{\E}{{\mathcal{E}}}

\renewcommand{\H}{{\mathcal{H}}}

 \newcommand{\M}{{\mathcal{M}}}
 
\renewcommand{\O}{{\mathcal{O}}}

 \newcommand{\T}{{\mathcal{T}}}

\newcommand{\upchi}{{\raise.35ex\hbox{$\chi$}}}
\newcommand{\fK}{{\mathfrak{K}}}

\newcommand{\foral}{\text{ for all }}

\newcommand{\qand}{\quad\text{and}\quad}

\newcommand{\qfor}{\quad\text{for}\quad}

\newcommand{\Aut}{\operatorname{Aut}}

\newcommand{\spn}{\operatorname{span}}

\newcommand{\rowt}{(T_e)_{e\in E^1}}

\newcommand{\rowl}{(L_e)_{e\in E^1}}

\begin{document}

\title[Limit Algebras and Directed Graphs]%
{A Class of Limit Algebras associated with Directed Graphs}
%
\author[D.W.Kribs and B. Solel]{David~W.~Kribs \, and \, Baruch~Solel}
\thanks{{\it key words and phrases.} directed graph, periodic weighted shift, Fock space,
limit algebra, Toeplitz-Cuntz-Krieger algebra, Cuntz-Pimsner
algebra, factor map.}
\address{Department of Mathematics and Statistics, University of Guelph,
Guelph, Ontario, Canada  N1G 2W1} \email{dkribs@uoguelph.ca}
\address{Department of Mathematics, Technion - Israel Institute of Technology, Haifa 32000, Israel}
\email{mabaruch@techunix.technion.ac.il}

\date{}
\begin{abstract}
Every directed graph defines a Hilbert space and a family of
weighted shifts that act on the space.  We identify a natural
notion of periodicity for such shifts and study their
$\ca$-algebras. We prove the algebras generated by all shifts of a
fixed period are of Cuntz-Krieger and Toeplitz-Cuntz-Krieger type.
The limit $\ca$-algebras determined by an increasing sequence of
positive integers, each dividing the next, are proved to be
isomorphic to Cuntz-Pimsner algebras and the linking maps are
shown to arise as factor maps. We derive a characterization of
simplicity and compute the $K$-groups for these algebras. We prove
a classification theorem for the class of algebras generated by
simple loop graphs.
\end{abstract}
\maketitle


\section{Introduction}\label{S:intro}

In this paper we initiate the study of a new class of
$\ca$-algebras associated with directed graphs. There is a family
of weighted shift operators associated with every directed graph
and, after identifying a natural notion of periodicity for these
shifts, we conduct an in-depth analysis of their associated
$\ca$-algebras. Specifically, we explicitly identify the structure
of the $\ca$-algebra generated by all shifts of a given period and
the limit algebras obtained from increasing sequences of positive
integers, each dividing the next, strictly in terms of familiar
objects from modern operator algebra theory.

Our initial motivation derives from work of Bunce and Deddens
\cite{BD1,BD2} from over thirty years ago in which a class of
$\ca$-algebras was studied via a limit algebra construction that
involved algebras generated by periodic weighted shift operators
on Hilbert space. The Bunce-Deddens algebras have proved to be an
extremely useful concrete class of operator algebras and have
arisen in a number of diverse settings
\cite{Archbold,Bre,DG,D,EGLP,Evans,Kbilateral,Kinductive,Moto,Pasnicu,PowerBD,Rordam1,SolelBD}.
We were also motivated by recent work of the first author
\cite{Kinductive}, where a generalization of this class was
recently obtained for the setting of Cuntz and Toeplitz-Cuntz
algebras. As we show, the class of algebras studied here contains
the Bunce-Deddens algebras and the algebras from \cite{Kinductive}
as the subclass generated by single vertex directed graphs with
$k$ loop edges (for $k=1$ and $k\geq 2$ respectively).

Our investigations draw on numerous aspects of contemporary
operator algebra theory. We make use of fundamental results from
the theory of graph $\ca$-algebras \cite{cbms} and Cuntz-Pimsner
algebras \cite{Pim,MS2,FMR,Ka2,MT}.  The theory of $\ca$-algebras
associated with `topological graphs', introduced by the second
author and Muhly \cite{MS_Morita} and studied further in
\cite{Ka1,Ka2,Ka3,MT} plays a central role. We utilize the theory
of `factor maps' recently invented by Katsura \cite{Ka2}. Each of
these tools complements our predominantly spatial analysis and,
further, our `Fock space' terminology indicates a philosophical
connection with theoretical physics.

The next section ($\S~\ref{S:prelim}$) contains requisite
preliminary material on graph $\ca$-algebras. We describe how
weighted shifts arise from directed graphs $E$ in
$\S~\ref{S:wtdshifts}$ and in $\S~\ref{S:periodicity}$ we identify
an appropriate notion of periodicity for these shifts, determined
by the path structure of the graph. The rest of the paper consists
of a detailed analysis of the $\ca$-algebras associated with
periodic shifts. In $\S~\ref{S:B(n)}$ we prove the algebras
$\A(n)$ and $\B(n)$ generated by shifts of a given period are of
Cuntz-Krieger and Toeplitz-Cuntz-Krieger type in such a way that
the explicit connection with the underlying graph is evident. Then
($\S~\ref{S:factor}$, $\S~\ref{S:limitalgebras}$) we identify the
corresponding limit algebras $\B_E(\{n_k\})$ as Cuntz-Pimsner
algebras $\O(E(\infty))$, where the topological graph $E(\infty)$
is defined by the path structure of $E$ and the sequence
$\{n_k\}$. In $\S$~\ref{S:example} we prove a classification
theorem for the algebras $\B_{C_j}(\{n_k\})$ generated by simple
loop graphs $C_j$. We compute the $K$-groups for the algebras
$\B_E(\{n_k\})$ in $\S$~\ref{S:Ktheory}. We finish in
$\S$~\ref{simp} by deriving a characterization of simplicity for
$\B_E(\{n_k\})$ in terms of $E(\infty)$ and discuss the connection
with $E$. Throughout the paper we have included a number of open
problems motivated by this work.

\section{Directed Graphs and Their $C^*$-Algebras}\label{S:prelim}

Let $E=(E^0,E^1,r,s)$ be a directed graph with vertices $x\in
E^0$, directed edges $e\in E^1$ and range and source maps $r,s$
giving the final and initial vertices of a given directed edge. We
shall assume $E$ is {\it finite and has no sources and no sinks},
so that every vertex in $E^0$ is the initial vertex for some  edge
and the final vertex for some edge. The finiteness assumption is
motivated by the $\ca$-algebra setting we will work in and the no
sink assumption is motivated by our definition of periodicity. We
focus on graphs with no sources simply to streamline the
presentation (see Remark~\ref{Etilderemark}).

Let $E^*$ be the set of all finite paths in $E$ and include the
vertices $E^0$ in $E^*$ as trivial paths. Given a path $w $ in $E$
we write $w = ywx$ when the initial and final vertices of $w$ are
$s(w)=x$ and $r(w)=y$ respectively. Thus, given $w_1,w_2\in E^*$
with $s(w_2)=r(w_1)$, $w_2 w_1$ is the path obtained by first
`walking' along $w_1$ and then along $w_2$. For $w$ in $E^*$ we
write $|w|$ for the length of $w$ and put $|x| = 0$ for every
vertex $x\in E^0$. Given $n \geq 0$, let $E^{=n}$ be the set of
paths in $E^*$ of length $n$, so that
\[
E^{=n} = \{ w\in E^* : |w| = n \}.
\]

There are two important $C^*$-algebras associated with every such
graph: the Cuntz-Krieger algebra $C^*(E)$ (or $\O(E)$) and its
Toeplitz extension $\T(E)$. For a recent survey of these algebras
we point the reader to the article \cite{cbms}. Both $\O(E)$ and
$\T(E)$ can be described either as universal objects or
concretely. We start by recalling their universal properties.

Given a directed graph $E$, a family $\{P_x,\, S_e : x\in E^0,\,
e\in E^1 \}$ of projections (one for each vertex) and partial
isometries (one for each edge) is said to be a
\emph{Toeplitz-Cuntz-Krieger $E$-family} (or a \emph{TCK
$E$-family} for short) if it satisfies the relations
\[
(\dagger)  \left\{
\begin{array}{lll}
(1)  & P_x P_y = 0 & \mbox{$\forall\, x,y \in E^0$, $ x \neq y$}
\\ (2) & S_{e}^{*}S_f = 0 & \mbox{$\forall\, e, f \in E^1$, $e
\neq f $}  \\ (3) & S_{e}^{*}S_e = P_{s(e)} & \mbox{$\forall\, e
\in E^1$}      \\ (4)  & \sum_{r(e)=x}\, S_e S_{e}^{*} \leq P_{x}
& \mbox{$\forall\, x \in E^0.$}
\end{array}
\right.
\]
Also, such a family is said to be a \emph{Cuntz-Krieger
$E$-family} (or a \emph{CK $E$-family}) if  equality holds in
$(4)$ whenever the set $r^{-1}(x)$ is non-empty.

The $C^*$-algebra $\O(E)$ is generated by a CK $E$-family
$\{p_x,s_e\}$ and has the property that, whenever $\{P_x,S_e\}$ is
a CK $E$-family inside a $C^*$-algebra $\B$, there is a
$*$-homomorphism $\pi_{P,S}$ from $\O(E)$ into $\B$ carrying $p_x$
to $P_x$ and $s_e$ into $S_e$. The Toeplitz algebra $\T(E)$ has a
similar universal property but with TCK $E$-families replacing CK
$E$-families.

It is also convenient to consider  concrete constructions of these
algebras. The details of the construction will be important when
we define the generalized Bunce-Deddens algebras through a spatial
approach.

Let $\H_E = \ell^2(E^*)$ be the Hilbert space with orthonormal
basis $\{ \xi_w : w\in E^* \}$ indexed by elements of $E^*$.
Define a family of partial isometries on $\H_E$ as follows: For
each $v\in E^*$ let
\begin{eqnarray}\label{ledefn}
L_v \, \xi_w = \left\{ \begin{array}{cl} \xi_{vw} & \mbox{if $s(v)
=
r(w)$} \\
0 & \mbox{if $s(v)\neq r(w)$}
\end{array}\right..
\end{eqnarray}
We will use the convention $\xi_{vw} = 0$ when $r(w)\neq s(v)$. We
shall put $ L_{x} \equiv P_x$ for the vertex projections.

 From a physical perspective the $L_e$ may be regarded as partial
`creation operators' acting on a generalized `Fock space' Hilbert
space $\H_E$. Alternatively, from a representation theory point of
view the $L_e$  arise from the left regular representation of
$E^*$ on $\H_E$. Evidently the family $\{P_x,L_e\}$ form a TCK
$E$-family. In fact, the $*$-homomorphism $\pi_{P,L}$ determined
by  the left regular representation is a $*$-isomorphism of
$\T(E)$ onto the $C^*$-algebra generated by the operators
$\{L_e\}$ \cite{FR,cbms}. Thus, for our purposes, we may identify
the algebra $\T(E)$ with this faithful representation
$\pi_{P,L}(\T(E))$. We shall, therefore,  for the sake of brevity
define the Toeplitz algebra (concretely) as follows.

\begin{defn}
The {\it Toeplitz algebra} of $E$ is the $\ca$-algebra $$\T(E)
\equiv \ca(\{L_w:w\in E^* \}) = \ca(\{L_e: e\in E^1\}).$$
\end{defn}

Let $R_v$, $v\in E^*$, be the partial isometries on $\H_E$
determined by the right regular representation of $E^*$, so that
$R_v \xi_w = \xi_{wv}$. It is easy to see that the subspaces $R_x
\H_E$ are invariant for $\T(E)$. (In fact, this family of
subspaces forms the unique maximal family of mutually orthogonal
irreducible subspaces for the associated `free semigroupoid
algebra' \cite{KP1}.)

\begin{prop}\label{compacts}
Let $\fK$ be the set of compact operators on $\H_E$. Then $\T(E)$
contains the subalgebra of compact operators $\fK_E =
\bigoplus_{x\in E^0} R_x \fK R_x$.
\end{prop}

\Prf By assumption $E^0$ is finite and hence $\T(E)$ is unital as
$I = \sum_{x\in E^0} P_x$. Note that for all $x\in E^0$, the rank
one projection $\xi_x \xi_x^*$ onto the subspace spanned by
$\xi_x$ satisfies
\[
\xi_x \xi_x^* = P_x \Big( \sum_{y\in E^0} \xi_y \xi_y^* \Big) =
P_x \Big( I - \sum_{e\in E^1} L_e L_e^* \Big).
\]
Thus, each $\xi_x \xi_x^*$ belongs to $\T(E)$. For an arbitrary
matrix unit $\xi_v \xi_w^*$ with $s(v)=x=s(w)$ we have
\[
\xi_v\xi_w^* = L_v \Big( \xi_{x} \xi_{x}^* \Big) L_w^* \in \T(E),
\]
and it follows that $\T(E)$ contains each $R_x \fK R_x$. \bx

Given a scalar $z\in\bbT$ we may define a {\it gauge unitary}
$U_z\in \B(\H_E)$ via
\begin{eqnarray}\label{gaugedefn}
U_z \xi_w = z^{|w|} \xi_w \qfor w\in E^*.
\end{eqnarray}
Then $\beta_z (L_e) = U_z L_e U_z^* = z L_e$ defines an
automorphism of $\T(E)$. Moreover, this automorphism leaves the
ideal $\fK_E$ invariant and hence factors through to an
automorphism on the quotient algebra $\T(E) / \fK_E$. It follows
that there is a continuous gauge action $\beta: \bbT \rightarrow
\Aut (\T(E) / \fK_E)$ and we obtain the following well-known
result.

\begin{thm}\label{gaugeinv}
The quotient algebra $\T(E) / \fK_E$ is isomorphic to the
universal Cuntz-Krieger algebra $\O(E)$.
\end{thm}

\Prf The operators $\{ L_e +\fK_E, P_x + \fK_E\}$ form a CK
$E$-family. Further, every vertex projection $P_x$ has infinite
rank since $E$ has no sources. Hence $P_x + \fK_E$ is a non-zero
projection for all $x\in E^0$. The result now follows from the
`gauge-invariant uniqueness theorem' for Cuntz-Krieger algebras
\cite{cbms}.
 \bx

\begin{rem}
Let us note that in the case of single vertex graphs, with $k\geq
1$ loop edges, the algebras $\T(E)$ and $\O(E)$ are familiar
objects. For $k=1$, $\T(E)$ is the $\ca$-algebra generated by the
unilateral shift $U_+$, or equivalently  the algebra of Toeplitz
operators with continuous symbol $\T_{C(\bbT)}$
\cite{Douglas_text}. For $k \geq 2$, $\T(E)$ is the Toeplitz-Cuntz
algebra $\E_k$ generated by the left creation operators $L_1,
\ldots , L_k$ which act on unrestricted $k$-variable Fock space.
On the other hand, $\O(E)$ is the algebra of continuous functions
on the unit circle $C(\bbT)$ ($k=1$) and the Cuntz algebra $\O_k$
($k\geq 2$). For more general graphs, analysis of $\T(E)$, $\O(E)$
and related non-selfadjoint algebras has occupied a substantial
portion of the recent operator algebra literature.
\end{rem}


\section{Weighted Shifts}\label{S:wtdshifts}

Consider a directed graph $E =(E^0,E^1,r,s)$.

\begin{defn}
A family of operators $\rowt$ that act on a Hilbert space $\H$ is
a {\it weighted shift} if there is a unitary $U: \H_E \rightarrow
\H$ and scalars $\Lambda = \{\lambda(w) : w\in E^*\setminus E^0\}$
such that the operators $(U^* T_e U)_{e\in E^1}$ satisfy
\begin{eqnarray}\label{shiftdefn}
(U^* T_e U) \xi_w = \lambda(ew) \xi_{ew} \quad\foral\, e\in E^1,
\, w\in E^*.
\end{eqnarray}
\end{defn}

We use the convention $\lambda(ew) =0 $ if $s(e)\neq r(w)$. To
streamline the discussion, we shall assume $\rowt$ acts on $\H_E$
and the unitary $U$ in the above definition is the identity
operator. The following pair of results are simple generalizations
of the single variable case.

\begin{prop}\label{basicprops}
Let $\rowt$ be a weighted shift with weights $\Lambda =
\{\lambda(w) :w\in E^*\setminus E^0\}$.
\begin{itemize}
\item[$(i)$] Then each $T_e$ factors as $T_e = L_e W_e$ where $L_E
= \rowl$ is the unweighted shift on $\H_E$ and each $W_e$ is a
positive operator which is diagonal with respect to the standard
basis for $\H_E$. \item[$(ii)$] There is a unitary $U\in
\B(\H_E)$, which is diagonal with respect to the standard basis
for $\H_E$, such that the weighted shift $(U^* T_e U)_{e\in E^1}$
has weights $ \{|\lambda(w)| :w\in E^*\setminus E^0\}$.
\end{itemize}
\end{prop}

\Prf For $(i)$, if we are given $e\in E^1$ then simply define $W_e
\xi_w = \lambda(ew) \xi_w$ if $s(e) = r(w)$ and $W_e \xi_w = 0$
otherwise. To see $(ii)$, define a unitary $U$ on $\H_E$ by
$U\xi_w = \mu(w)\xi_w$ where $\mu(x) = 1$ for $x\in E^0$ and for
$|w|>1$, say $w=ev$, $|\mu(w)| =1$ is chosen so that
$(\ol{\mu(v)}\,\lambda(w))\,\mu(w) \geq 0.$ Then if $s(e) = r(v)$
we have
\begin{eqnarray*}
\big( UT_eU^* \big) \, \xi_v &=& \ol{\mu(v)} \, UT_e  \xi_v =
\ol{\mu(v)}\lambda(w) U \xi_{ev} \\ &=& \big( \ol{\mu(v)}
\lambda(w) \mu(w) \big) \, \xi_{ev} = |\lambda_w| \, \xi_{ev}.
\end{eqnarray*}
Whereas $UT_eU^* \xi_v =0$ if $s(e)\neq r(v)$.  \bx

\section{Periodicity}\label{S:periodicity}

Now let $E =(E^0,E^1,r,s)$ be a fixed finite graph with no sinks
and no sources and let $n \geq 1$ be a fixed positive integer.

Observe that every $w\in E^*$ has a unique factorization of the
form
\[
w = w(n) \,\, v_k \cdots v_1
\]
where $v_i \in E^{=n}$  for $1\leq i \leq k$ and $w(n)\in E^{<
n}$.

\begin{defn}\label{periodicshiftdefn}
A weighted shift $T= \rowt$ with weights $\Lambda = \{ \lambda(w)
: w\in E^*\setminus E^0\}$ is {\it period $n$} if $\lambda(w) =
\lambda(w(n))$ for all $\lambda(w) \in \Lambda$. In other words,
\begin{eqnarray*}
T_e \xi_w = \lambda(ew) \xi_{ew} = \lambda((ew)(n)) \xi_{ew} \quad
\forall \,e\in E^1, \, w\in E^*.
\end{eqnarray*}
\end{defn}

The weights $\lambda(w) \in \Lambda$ with $|w| \leq n$ are the
`remainder weights' for an $n$-periodic shift.

\begin{defn}
Let $\A(n)$ be the $\ca$-algebra generated by the $T_e$, $e\in
E^1$, from all $n$-periodic weighted shifts $T = (T_e)_{e\in E^1}$
on $\H_E$. Let $\{ n_k \}_{k\geq 1}$ be an increasing sequence of
positive integers such that $n_k | n_{k+1}$ for $k \geq 1$.
Observe that every period $n_k$ weighted shift $T=\rowt$ is also
period $n_{k+1}$. Thus,
\[
\A(n_{1}) \subseteq \A(n_2) \subseteq \ldots \subseteq \A(n_k)
\subseteq \ldots
\]
and we may consider the (norm-closed) limit algebra
\[
\A_E (\{n_k\}) := \ol{\bigcup_{k\geq 1} \A(n_k)}.
\]
As $\A(n)$ contains the $\ca$-algebra $\T_E$ generated by the
unweighted shifts $L_E = (L_e)_{e\in E^1}$, by
Proposition~\ref{compacts} it contains the compact operators
$\fK_E$. Let $\B(n)$ be the quotient of $\A(n)$ by $\fK_E$, so
there is a short exact sequence
\[
0 \rightarrow \fK_E \rightarrow \A(n) \rightarrow \B(n)
\rightarrow 0.
\]
Thus, given a sequence $\{ n_k \}_{k\geq 1}$ we have the sequence
of injective inclusions
\[
\B(n_{1}) \subseteq \B(n_2) \subseteq \ldots \subseteq \B(n_k)
\subseteq \ldots,
\]
and we may also consider the limit algebra
\[
\B_E (\{n_k\}) := \ol{\bigcup_{k\geq 1} \B(n_k)}.
\]
We refer to $\B_E(\{n_k\})$ as a {\it generalized Bunce-Deddens
algebra} (or a {\it Bunce-Deddens graph algebra}).
\end{defn}

\section{The Algebras $\A(n)$ and $\B(n)$}\label{S:B(n)}

For the analysis of this section we shall fix a finite  graph $E$
with no sinks and no sources and a positive integer $n \geq 1$.
Recall that $E^{= n}$ is the set of paths in $E$ of length equal
to $n$. Similarly define $E^{< n}$ to be the set of paths in $E$
of length strictly less than $n$. Also define
\[
E^{{\rm per}\,n} = \{ w\in E^* : |w| = mn \,\,{\rm for
\,\,some}\,\, m\geq 0\}.
\]

\begin{defn}
Let $E(n) = (E(n)^0,E(n)^1,r_{E(n)},s_{E(n)})$ be the graph
defined as follows.

First define $E(n)^0 \equiv E^{<n}$. In other words, the paths of
length less than $n$ in $E$ now serve as the vertices of $E(n)$.
We may use $w$ to denote such a path in $E$ or a vertex in $E(n)$.
It will be clear from the context what is the role of $w$.
Moreover, when we write $r_E(w)$ we refer to the vertex in $E$
that is the range of the path $w$. This vertex can be viewed
either as a vertex of $E$ or as a vertex of $E(n)$ (since
$E^0\subseteq E(n)^0$) depending on the context.

Now we set
\[
E(n)^1=\{(e,w) \in E^1 \times E^{<n} : s_E(e)=r_E(w) \}
\]
and the maps $r_{E(n)}$ and $s_{E(n)}$ are defined by setting
\[
s_{E(n)}(e,w)=w
\] and
\[
r_{E(n)}(e,w)=\left\{ \begin{array}{cl} ew & \mbox{if $|w|<n-1$}
\\ r_E(e) & \mbox{if $|w|=n-1$} \end{array} \right.
\]
\end{defn}


We next define a TCK and a  CK $E(n)$-family. For this we first
let $T_{(e,w)}$, for $(e,w)\in E(n)^1$,  be the operator on
$\mathcal{H}_E$ defined by
\[
T_{(e,w)}\xi_{w'} = \left\{ \begin{array}{cl} \xi_{ew'} & \mbox{if
$w=w'(n)$}
\\ 0 & \mbox{if $w\neq w'(n)$} \end{array} \right..
\]
and $Q_w$, for $w\in E(n)^0$, be the projection onto the subspace
of $\mathcal{H}_E$ spanned by the vectors $\xi_{w'}$ with
$w'(n)=w$; so that
\[
Q_w = \sum_{w'(n)=w} \xi_{w'}\xi_{w'}^* = \sum_{v\in E^{{\rm
per}\,n},\, r(v)=s(w)} \xi_{wv}\xi_{wv}^*.
\]

Observe that
\[ T_{(e,w)}^*T_{(e,w)}=Q_w .
\]
It is also easy to check that
\[
T_{(e,w)}^*\xi_{w''}= \left\{ \begin{array}{cl} \xi_{v} & \mbox{if
$w''=ev,\, v(n)=w$}
\\ 0 & \mbox{otherwise} \end{array} \right..
\]

Thus, $ T_{(e,w)}T_{(e,w)}^*$ is the projection onto the subspace
spanned by all $\xi_{w''}$ with $w''=ev$ for $v$ satisfying
$v(n)=w$. It follows that, for $w_0\in E(n)^0$,
\[
\sum_{r(e,w)=w_0} T_{(e,w)}T_{(e,w)}^*= \left\{ \begin{array}{cl}
Q_{w_0} & \mbox{if $n>|w_0|>0$}
\\ Q_{w_0}- \xi_{w_0}\xi_{w_0}^* & \mbox{if  $|w_0|=0$} \end{array}
\right..
\]
Note the index set in this sum is a singleton whenever
$n>|w_0|>0$.

It follows that $\{Q_w,T_{(e,w)}\}$ is a TCK $E(n)$-family and,
thus, there is a $*$-homomorphism $\rho$ from $\T(E(n))$ into
$\B(\H_E)$ carrying the generators of the Toeplitz algebra to this
family.

Observe that every operator $T_{(e,w)}$ as above is the periodic
weighted shift associated to the weights $\Lambda_{(e,w)} = \{
\lambda(w') : w'\in E^* \setminus E^0\}$ where
\[
\lambda(w')=\left\{ \begin{array}{cl} 1 & \mbox{if $w'=ev,\, v(n)=w$} \\
0 & \mbox{otherwise} \end{array} \right. .
\]
Note also that the operators $T_e$ associated with every
$n$-periodic weighted shift $T=(T_e)_{e\in E^1}$ can be written as
a finite sum
\[
T_e=\sum_{w\in E^{<n},\, s_E(e)=r_E(w)} \lambda(ew)T_{(e,w)},
\]  It follows that the $*$-homomorphism $\rho$ described above
maps $\T(E(n))$ onto $\A(n)$. (This map is not an isomorphism
however, see the discussion following Theorem~\ref{Bn}).

Setting $S_{(e,w)}=q(T_{(e,w)})$ and $P_w=q(Q_w)$, where $q$ is
the quotient map from $\B(\mathcal{H}_E)$ onto the Calkin algebra,
we get a CK $E(n)$-family. Such a family defines a
$*$-homomorphism $\pi$ from $\mathcal{O}(E(n))$ into the Calkin
algebra. Since the $T_{(e,w)}$ generate $\A(n)$, the image of
$\pi$ is $\mathcal{B}(n)$. Of course, in principle the $P_w$ could
be zero.

\begin{lem}
For all $w\in E(n)^0$, $Q_w$ is an infinite rank projection, and
hence $P_w\neq 0$.
\end{lem}

\Prf Fix $w\in E(n)^0$. Since $E$ has no sources, we can find
paths $v_k \in E^{=nk}$ for $k\geq 1$ such that $s(w) = r(v_k)$.
Then $(wv_k)(n) = w$ for all $k \geq 1$. Thus, $\xi_{wv_k}$
belongs to the range of $Q_w$ for all $k\geq 1$.
 \bx

We may now prove the following.

\begin{thm}\label{Bn}
$\pi$ is a $*$-isomorphism of $\mathcal{O}(E(n))$ onto
$\mathcal{B}(n)$.
\end{thm}

\Prf  For $z \in \mathbb{T}$ let $U_z$ be the unitary operator on
$\mathcal{H}_E$ defined as in (\ref{gaugedefn}). Setting
$\gamma_z(R)=U_z R U_z^*$, we get a one-parameter semigroup of
(inner) automorphisms of $\B(\mathcal{H}_E)$. For $(e,w) \in
E(n)^1$, $w'\in E(n)^0$ and $z\in \mathbb{T}$,
\[
\gamma_z(T_{(e,w)})\xi_{w'}=U_z T_{(e,w)} U_z^* \xi_{w'}= z \,
T_{(e,w)}\xi_{w'} .
\]
Hence $\gamma_z(T_{(e,w)})=z\,T_{(e,w)}$ and it follows that each
$\gamma_z$ defines an automorphism of $\A(n)$. Moreover, as
discussed above, $\gamma_z$ leaves $\fK_E$ invariant and so
$\{\gamma_z \}$ induces a one parameter semigroup of automorphisms
on the quotient $\mathcal{B}(n)$, which we shall also denote by
$\{\gamma_z\}$.

Thus we have $\gamma_z(S_{(e,w)})=z\, S_{(e,w)}$ for all $(e,w)\in
E(n)^1$ and $\gamma_z(P_w)=P_w$ for all $w\in E(n)^0$. Since
$P_w\neq 0$ for all $w\in E(n)^0$, we can now apply the
gauge-invariant uniqueness theorem \cite{cbms} to conclude that
$\pi$ is an isomorphism.
 \bx

It is natural to ask whether there is a similar isomorphism result
for $\A(n)$ in place of $\B(n)$. One may try to use the
$*$-homomorphism $\rho: \T(E(n))\rightarrow \A(n)$ mentioned above
(defined via the TCK $E(n)$-family
 $\{Q_w,T_{(e,w)}\}$) but this map fails to be injective (for
 $n>1$) since, whenever $0<|w_0|<n$,
 $$\rho\Big(p_{w_0}-\sum _{r(e,w)=w_0}
 s_{(e,w)}s_{(e,w)}^*\Big)=Q_{w_0}-\sum_{r(e,w)=w_0} T_{e,w}T_{e,w}^*
 =0, $$ while $p_{w_0}-\sum _{r(e,w)=w_0}
 s_{(e,w)}s_{(e,w)}^*\neq 0$.

 This, of course, does not imply that $\T(E(n))$ is not
 $*$-isomorphic to $\A(n)$ via a different map. However, we
 will show in Note~\ref{noniso} that these algebras may be
 non-isomorphic.

Nevertheless, there is a way to get a result for $\A(n)$ that is
analogous to Theorem~\ref{Bn}. The inspiration for the analysis
sketched below comes from $\S$~3 of \cite{Ka2} and $\S$~7 of
\cite{MT}. In the terminology of \cite{MT}, the algebra $\A(n)$ is
a relative quiver algebra (see also Example 1.5 of \cite{FMR}),
related to the relative Cuntz-Pimsner algebras introduced in
\cite{MS2}.
 Since our main focus in this paper is on the
algebras $\B(n)$, and their direct limits, we shall only sketch
the construction and the results and leave some details to the
reader. The idea is to replace the graph $E(n)$ by another graph,
written $E[n]$. (Using the notation of \cite{Ka2}, $E[n]$ is
$E(n)_{E^0}$.)  To define it, we first let $c(E^0)$ be a copy of
$E^0$ (whose elements will be written $c(v)$, $v\in E^0$). Then
\[ E[n]^0=E(n)^0 \sqcup c(E^0) \] and \[
E[n]^1=E(n)^1 \sqcup \big\{(e,c(v)):e\in E^1, v\in E^0, s_E(e)=v
\big\}.\] The maps $s_{E[n]}$ and $r_{E[n]}$ coincide with
$s_{E(n)}$ and $r_{E(n)}$, respectively, on $E(n)^1$ and
\[s_{E[n]}(e,c(v))=c(v)\in c(E^0) \]
\[
r_{E[n]}(e,c(v))=\left\{ \begin{array}{cl} e\in E(n)^0 & \mbox{if
$n > 1$} \\ r_E(e) & \mbox{if $n=1$} \end{array}\right..
\]
The TCK $E(n)$-family $\{Q_w,T_{(e,w)}\}$ defined above gives rise
to a CK $E[n]$-family $\{G_u,R_z : u\in E[n]^0, z\in E[n]^1 \}$
defined by
\[
G_u = \left\{ \begin{array}{ll} Q_u-\xi_u\xi^*_u & \mbox{if } u
\in E^0 \\ \xi_v \xi_v^* & \mbox{if } u=c(v)\in c(E^0) \\ Q_u &
\mbox{if } u \in E(n)^0\backslash  E^0  \end{array} \right.
\] and \[
R_z=\left\{ \begin{array}{ll} T_{(e,v)}(Q_v-\xi_v\xi_v^*) &
\mbox{if } z=(e,v)\in E(n)^1, v\in E^0 \\ T_{(e,v)}\xi_v\xi_v^* &
\mbox{if } z=(e,c(v)), e\in E^1, v\in E^0 \\ T_{(e,w)} & \mbox{if
} z=(e,w), (e,w)\in E(n)^1 \backslash E^1 \end{array} \right. .\]

Each of the projections $G_u$, $u\in E[n]^0$, is non-zero and the
unitaries $U_z$, $z\in \bbT$, from (\ref{gaugedefn}) define a
semigroup of gauge automorphisms on the $\ca$-algebra generated by
$\{G_u,R_z\}$. Thus we may proceed as in Theorem~\ref{Bn} to show
the following.

\begin{thm}\label{A(n)}
The algebra $\A(n)$ is $*$-isomorphic to $\O(E[n])$.
\end{thm}

\vspace{5mm}

Even though the results above (Theorem~\ref{A(n)} and
Theorem~\ref{Bn}) provide descriptions of the algebras $\A(n)$ and
$\B(n)$, it is sometimes useful to have a more concrete
description. We shall now provide such descriptions in
Theorem~\ref{An} and Corollary~\ref{Bndecomp} below.

 Towards this end, let us define a directed graph
$E(\mbox{\footnotesize{=n}}) = (E(\mbox{\footnotesize{=n}})^0,$ $
E(\mbox{\footnotesize{=n}})^1,
r_{E(\mbox{\footnotesize{=n}})},s_{E(\mbox{\footnotesize{=n}})})$
as follows:
\[
E(\mbox{\footnotesize{=n}})^0 = E^0,
\]
\[
E(\mbox{\footnotesize{=n}})^1 = \big\{ (w,x): x\in E^0, \, w\in
E^{=n}, \, s(w)=x\big\},
\]
\[
s_{E(\mbox{\footnotesize{=n}})}((w,x)) = x, \quad
r_{E(\mbox{\footnotesize{=n}})} ((w,x)) = r_E(w).
\]
So the vertices of $E(\mbox{\footnotesize{=n}})$ are identified
with the vertices of $E^0$ and every path of length $n$ in $E$
determines an edge in $E(\mbox{\footnotesize{=n}})$.

For all $x\in E^0$ let $d_E(n,x)$ be the number of paths of length
strictly less than $n$ in $E$ with $x$ as initial vertex, and let
$d_E(n) = \sum_{x\in E^0} d_E(n,x)$ be the cardinality of the set
$E^{<n}$. Further let $\H_E (w)$, $w\in E^{<n}$, be the ranges of
the projections $Q_w$ discussed above, and so
\[
\H_E(w) = \spn \{ \xi_{wv} : s(w) = r(v),\,v\in E^{\rm{per}\,n}
\}.
\]
We shall write the set of paths for $E(\mbox{\footnotesize{=n}})$
as
\[
E(\mbox{\footnotesize{=n}})^* = \{(v,s(v)) : v\in E^{\rm{per}\,n}
\}
\]
and extend the range and source maps as
$s_{E(\mbox{\footnotesize{=n}})}((v,s(v))) = s_E (v)$ and
$r_{E(\mbox{\footnotesize{=n}})}((v,s(v))) = r_E (v)$. Note that
there is a natural correspondence between elements of
$E(\mbox{\footnotesize{=n}})^*$ and $E^{\rm{per}\,n}$.

The following theorem generalizes \cite[Theorem 4.1]{Kinductive}.
Let us write $P_x$, $x \in E(\mbox{\footnotesize{=n}})^0 = E^0$,
for the vertex projections of $\T(E(\mbox{\footnotesize{=n}}))$ in
the result below.

\begin{thm}\label{An}
$\A(n)$ is unitarily equivalent to the algebra
\begin{eqnarray}\label{Andecomp}
\quad\Big( \bigoplus_{x\in E^0} P_x^{(d_E(n,x))} \Big) \M_{d_E(n)}
\big( \T(E(\mbox{\footnotesize{=n}})) \big) \Big( \bigoplus_{x\in
E^0} P_x^{(d_E(n,x))} \Big),
\end{eqnarray}
where $\M_{d_E(n)} \big( \T(E(\mbox{\footnotesize{=n}})) \big)$ is
the full matrix algebra with entries in
$\T(E(\mbox{\footnotesize{=n}}))$ that acts on the space
$\bigoplus_{x\in E^0}
(P_x\H_{E(\mbox{\footnotesize{=n}})})^{(d_E(n,x))}$.
\end{thm}

\Prf For all $w\in E^{<n}$ define a unitary
\[
\begin{array}{rcl}
U_w : \H_E(w) & \,\longrightarrow\, & P_{s(w)}\H_{E(\mbox{\footnotesize{=n}})} \\
\xi_{wv} & \,\mapsto\, & \xi_{(v,s(v))}
\end{array}
\]
Since every $u\in E^*$ has a unique factorization of the form
$u=wv$ with $w\in E^{<n}$ and $v\in E^{\rm{per}\,n}$, it follows
that $\H_E$ decomposes as
\[
\H_E = \bigoplus_{w\in E^{<n}} \H_E(w).
\]
Thus we may define a unitary
\[
U: \H_E \longrightarrow \bigoplus_{x\in E^0} \Big( P_x
\H_{E(\mbox{\footnotesize{=n}})} \Big)^{(d_E(n,x))}
\]
by $ U = \bigoplus_{w\in E^{<n}} U_w. $

Let $T_{(e,v)}$ be a generator of $\A(n)$, as discussed prior to
Theorem~\ref{Bn}, so that $v\in E^{<n}$ and $s(e) =r(v)$. Then for
$w,w'\in E^{<n}$ we have
\[
U_{w'} T_{(e,v)} U_w^* = \left\{ \begin{array}{cl} 0 & \mbox{if
$v\neq w$ or $ev\neq w'$} \\ P_{s(v)} & \mbox{if $v= w$, $w'=ev$
and $0\leq |v| < n-1$} \\ L_{(w',s(v))} & \mbox{if $v= w$, $w'=ev$
and $|v| = n-1$} \end{array}\right.
\]
where here $P_{s(v)}$ and $L_{(w',s(v))}$ are generators for
$\T(E(\mbox{\footnotesize{=n}}))$. It follows that $U\A(n) U^*$ is
given by the matrix algebra in (\ref{Andecomp}).
 \bx

\begin{note}\label{noniso} As we noted above,
the homomorphism $\rho : \T(E(n))\rightarrow \A(n)$ which induces
the $*$-isomorphism $\pi$ above is not an isomorphism. In fact, if
we consider the case that $E$ is the single vertex, single loop
edge graph and $n=2$, we can see that $\T(E(2))$ and $\A(2)$ are
different algebras. To see this we argue as follows. First note
that, from the previous theorem, it follows that $\A(2)$ is
isomorphic to $\M_2(\T_{C(\bbT)})$. Note that every ideal of
$\M_2(\T_{C(\bbT)})$ is of the form $\M_2(J_0)$ for some ideal
$J_0$ in $\T_{C(\bbT)}$. But then every (non-trivial) quotient of
 $\M_2(\T_{C(\bbT)})$ is of the form $\M_2(\A)$ for some
 commutative $C^*$-algebra $\A$. As a result, every irreducible
 representation of every (non-trivial) quotient of $\A(2)$ is
 finite-dimensional. This is not true for the algebra $\T(E(2))$.
 (Simply consider the quotient of this algebra by the kernel of
 the map $\rho$ discussed above. The quotient is isomorphic to
 $\M_2(\T_{C(\bbT)})$ and this algebra can be represented
 irreducibly on the orthogonal sum of two copies of $H^2$.)
 It follows that $\T(E(2))$ is not isomorphic to $\A(2)$.

\end{note}

Letting $p_x$, $x\in E(\mbox{\footnotesize{=n}})^0$, denote the
projections from a universal family $\{s_e,p_x\}$ for
$\O(E(\mbox{\footnotesize{=n}}))$, we obtain the following
alternative characterization of $\B(n)$.

\begin{cor}
There is a $*$-isomorphism from $\B(n)$ onto the algebra
\begin{eqnarray}\label{Bndecomp}
\quad \Big( \bigoplus_{x\in E^0} p_x^{(d_E(n,x))} \Big)
\M_{d_E(n)} \big( \O(E(\mbox{\footnotesize{=n}})) \big) \Big(
\bigoplus_{x\in E^0} p_x^{(d_E(n,x))} \Big).
\end{eqnarray}
\end{cor}

\begin{rem}
This result generalizes the fact that in the case of single vertex
graphs with $k\geq 1$ loop edges, $\B(n)$ is isomorphic to $\M_n
(C(\bbT))$ for $k=1$ \cite{BD1,BD2} and $\M_{\frac{k^n -1}{k-1}}
(\O_{k^n})$ for $k\geq 1$ \cite{Kinductive}. Moreover, by
Theorem~\ref{Bn} this result yields what is, to our knowledge, a
new description of universal graph $\ca$-algebras of the form
$\O(E(n))$.
\end{rem}

\begin{cor}
Given a directed graph $E$ and a positive integer $n\geq 1$, there
is a $*$-isomorphism from $\O(E(n))$ onto the $\ca$-algebra in
(\ref{Bndecomp}).
\end{cor}

\section{Factor Maps}\label{S:factor}

Let us examine in more detail the embedding maps that determine
the limit algebras $\B_E(\{n_k\})$. Fix $n,k \in \mathbb{N}$ and
write $\pi_n : \mathcal{O}(E(n)) \rightarrow \mathcal{B}(n)$ and
$\pi_{nk} : \mathcal{O}(E(nk)) \rightarrow \mathcal{B}(nk)$ be the
$*$-isomorphisms of Theorem~\ref{Bn}.

Recall that the algebra $\mathcal{A}(n)$ is contained in
$\mathcal{A}(nk)$. We write $i_{nk,n}$ (or, simply, $i$) for this
inclusion map and $\bar{i}_{nk,n}$ for the embedding
$\bar{i}_{nk,n} : \mathcal{B}(n) \rightarrow \mathcal{B}(nk)$
induced by $i_{nk,n}$. Letting $j_{nk,n}=\pi_{nk}^{-1} \circ
\bar{i}_{nk,n} \circ \pi_n$ we get an injective $*$-homomorphism
\begin{equation}\label{j}
j_{nk,n} : \O(E(n)) \rightarrow \O(E(nk)).
\end{equation}

For $(e,w) \in E(n)^1$, $T_{(e,w)}$ is an operator from a shift of
period $n$ and
\[
i_{nk,n}(T_{(e,w)})=\sum_{(e,w') \in E(nk)^1,\,
w'(n)=w}T_{(e,w')}.
\]
Thus
\begin{eqnarray}\label{jeqn}
j_{nk,n}(S_{(e,w)})=\sum_{(e,w') \in E(nk)^1,\,
w'(n)=w}S_{(e,w')},
\end{eqnarray}
where here $S_{(e,w)}$ and $S_{(e,w')}$ are generators of
$\O(E(n))$ and $\O(E(nk))$ respectively.

We now show that the map $j_{(nk,n)}$ is induced from a `regular
factor map' $m: E(nk) \rightarrow E(n)$. Factor maps were
introduced and studied recently by Katsura  \cite{Ka2} in the
context of topological graphs whose vertex and edge spaces are
locally compact topological spaces. Here we shall need these
concepts only for finite graphs and, thus, the definitions can be
simplified.

\begin{defn}\label{factor}
Let $F=(F^0,F^1,s_F,r_F)$ and $E=(E^0,E^1,s_E,r_E)$ be finite
graphs. A {\it factor map} from $F$ to $E$ is a pair $m=(m^0,m^1)$
consisting of maps $m^0:F^0 \rightarrow E^0$ and $m^1:F^1
\rightarrow E^1$ such that
\begin{enumerate}
\item[$(i)$] For every $e\in F^1$,
\[
r_E(m^1(e))=m^0(r_F(e)) \qand s_E(m^1(e))=m^0(s_F(e)).
\]
\item[$(ii)$] If $e'\in E^1$ and $v\in F^0$ satisfy
$s_E(e')=m^0(v)$, then there exists a unique element $e\in F^1$
such that $m^1(e)=e'$ and $s_F(e)=v$.
\end{enumerate}
Such a map is said to be {\it regular} if also
\begin{enumerate}
\item[$(iii)$] $(r_F)^{-1}(v)$ is non-empty whenever $v\in F^0$
and $(r_E)^{-1}(m^0(v))$ is non-empty.
\end{enumerate}
\end{defn}

We now define
\[ m=(m^0,m^1): E(nk) \rightarrow E(n) \] by
\[ m^0(w)=w(n)\;,\;\; w\in E(nk)^0 \] and
\[ m^1(e,w)=(e,m^0(w))\;,\;\; (e,w)\in E(nk)^1. \]

\begin{lem}\label{m}
The pair $m=(m^0,m^1)$ defined above is a regular factor map from
$E(nk)$ to $E(n)$.
\end{lem}

\Prf First note that, for $(e,w)\in E(nk)^1$, $w(n)$ is indeed in
$E(n)^0$ (as $|w(n)|<n$) and $(e,w(n))$ is in $E(n)^1$ (as
$s_E(e)=r_E(w)=r_E(w(n))$).

Fix $(e,w)\in E(nk)$. Then
\[
s_{E(n)}(m^1(e,w))=s_{E(n)}(e,w(n))=w(n)=m^0(w)= s_{E(nk)}(e,w).
\]
To prove a similar statement for $r$ in place of $s$ we
distinguish two cases: when $|w|<nk-1$ (and so $|w(n)|<n-1$) and
when $|w|=nk-1$ (and $|w(n)|=n-1$). In the first case
\begin{eqnarray*}
r_{E(n)}(m^1(e,w))&=& r_{E(n)}(e,w(n))\\ &=&(ew)(n) = m^0(ew)=
m^0(r_{E(nk)}(e,w))
\end{eqnarray*}
and in the latter
\begin{eqnarray*}
r_{E(n)}(m^1(e,w))&=& r_{E(n)}(e,w(n)) \\ &=& r_{E}(e) =
m^0(r_{E}(e))=m^0( r_{E(nk)}(e,w)).
\end{eqnarray*}
This establishes part (i) of the definition.

For (ii), suppose $(e,w')\in E(n)^1$ and $w_1 \in E(nk)^0$ satisfy
$w' = s_{E(n)}(e,w')=m^0(w_1)=w_1(n)$. Then $(e,w_1)$ lies in
$E(nk)^1$, as $r_E(w_1)=r_E(w_1(n))=r_E(w')=s_E(e)$, and satisfies
$m^1(e,w_1)=(e,w_1(n))=(e,w')$ and $s_{E(nk)}(e,w_1)=w_1$, proving
part (ii).

The claim that the map is regular follows from the fact that
$E(nk)$ has no sources, since $E$ has none.
 \bx

The following proposition is \cite[Proposition 2.9]{Ka2} applied
to finite graphs.

\begin{prop}\label{reg}
Let $E$ and $F$ be two finite graphs and $m$ be a regular factor
map from $F$ to $E$. Then there is a unique $*$-homomorphism
$\mu_m :\O(E) \rightarrow \O(F)$ such that, for every $v\in E^0$
and $e \in E^1$,
\begin{enumerate}
\item[(1)]
 $\mu_m(P_v)=\sum_{u\in (m^0)^{-1}(v)} P_{u}$
 \item[and]
\item[(2)] $\mu_m(S_e)=\sum_{f\in (m^1)^{-1}(e)}S_f$.
\end{enumerate}
Also, $\mu_m$ is injective if and only if $m^0$ is surjective.
\end{prop}

Returning to $E(nk)$ and $E(n)$, together with (\ref{jeqn}) the
above immediately implies the following.

\begin{cor}\label{jm}
The regular factor map $m$ of Lemma~\ref{m} satisfies $$
j_{nk,n}=\mu_m. $$
\end{cor}

\begin{rem}\label{factorq}
Replacing the graphs $E(nk)$ and $E(n)$ by $E[nk]$ and $E[n]$
respectively (as in the discussion leading to Theorem~\ref{A(n)}),
one can define a factor map $q=(q^0,q^1)$ from $E[nk]$ to $E[n]$
where $q^i$ agrees with $m^i$ on $E(nk)^i$, $i=1,2$,
$q^0(c(v))=c(v)$ and $q^1(e,c(v))=(e,c(v))$. As in
Corollary~\ref{jm}, the map $\mu_q$ induced by $q$ is the
embedding of $\O(E[n])$ into $\O(E[nk])$ induced by the embedding
of $\A(n)$ into $\A(nk)$.
\end{rem}

\section{$\B_E(\{n_k\})$ as a Cuntz-Pimsner Algebra}\label{S:limitalgebras}

Fix a finite graph $E=(E^0,E^1,r_E,s_E)$ with no sinks and no
sources and an increasing sequence $\{n_k \}_{k\geq 1}$ of
 positive integers with
each $n_k$ dividing $n_{k+1}$ (and write $m_k=n_{k+1}/n_k$). We
also write $n_0=1$.

It will be important for us to note that every $w\in E^*$ with
$|w|=m$ can be written uniquely as
\begin{eqnarray} \label{decomp}
w=w_1 w_2 \cdots w_k
\end{eqnarray}
where
\begin{eqnarray}\label{xi}
w_i \in X_i \equiv \big\{ w\in E^* :
 0 \leq |w|< n_i,\, |w|\equiv 0\; (mod \;n_{i-1}) \big\}.
\end{eqnarray}
Supposing that  $|w_i| = k_i n_{i-1}$, and so $k_i < m_{i-1}$, we
have
\begin{eqnarray}\label{wi}
|w_i|=k_{i}n_{i-1}\leq n_i-n_{i-1}
\end{eqnarray}
and
\begin{eqnarray*}
m=\sum_{i=1}^k k_{i}n_{i-1}.
\end{eqnarray*}
Note that (\ref{wi}) holds for all $w\in X_i$.

Now write $X$ for the (compact) product space $X=X_1\times X_2
\times \cdots $ and $Y$ for the (closed) subset
\[ Y=\big\{\omega =(w_1,w_2, \ldots ) \in X: s_E(w_k)=r_E(w_{k+1}),\;k=1,2
\cdots \big\} .\] Also, let
 $\tau : E^* \rightarrow Y$
be the map defined by \begin{eqnarray} \label{tau}
\tau(w)=(w_1,w_2, \ldots ,w_k, s_E(w_k),s_E(w_k), \ldots ) \in Y,
\end{eqnarray}
where  $w=w_1 w_2 \cdots w_k$  is the decomposition as in
(\ref{decomp}) and (\ref{xi}). Then $\tau$ is an embedding of
$E^*$ onto a dense subset of $Y$. We shall refer to $Y$ as
\emph{the $\{n_k\}$-compactification of $E^*$}.

\begin{defn}
For $e\in E^1$, we define {\it odometer maps} $\sigma_e: D_e
\rightarrow R_e$ on $Y$ as follows. First, put
\[ D_e =\big\{y=(y_1,y_2,\ldots )\in Y: r_E(y_1)=s_E(e) \big\} \]
and
\[ R_e =\big\{y=(y_1,y_2,\ldots )\in Y : \mbox{for some}\;l\leq \infty, \;
 y_i=r_E(e)
\;\mbox{for all}\; i<l \] \[  \mbox{and (if}\;l\neq \infty )\;
y_l=ew' \; \mbox{for some}\; |w'|\equiv -1 (mod\; n_{l-1}) \big\}.
\]
Now, given $\omega =(w_1,w_2,\ldots )$ in $D_e$, and recalling
(\ref{wi}), write $i(w)$ for the smallest positive integer $i$
such that $|w_i|<n_i-n_{i-1}$ (if there is one) or $i(w)=\infty$
if $|w_i|=n_i-n_{i-1}$ for every $i$.

If $i(w)<\infty$ we write $\sigma_e(w)=u$ where $u_i=r_E(e)$ for
$i<i(w)$, $u_{i(w)}=ew_1w_2 \cdots w_{i(w)}$ and $u_i=w_i$ for
$i>i(w)$. If $i(w)=\infty$, we set $\sigma_e(w)=(r_E(e),r_E(e),
\ldots )$.
\end{defn}

\begin{lem}\label{sigma}
For every $e\in E^1$,
\begin{enumerate}
\item[(i)] $\tau(ew)=\sigma_e(\tau(w))$ for every $w\in E^*$ such
that $s_E(e)=r_E(w)$. \item[(ii)] The sets $D_e$ and $R_e$ are
compact and $\sigma_e$ is a continuous map from $D_e$ onto $R_e$.
\end{enumerate}
\end{lem}

\Prf Let $e\in E^1$. To prove (i), fix $w\in E^*$ with $s_E(w) =
r_E(w)$ and let $m = i(w)$. Then $w=w_1w_2 \cdots w_k$ (as in
(\ref{decomp}) and (\ref{xi})) and $|w_i|=n_i-n_{i-1}$ for $i<m$
and $|w_m|<n_m-n_{m-1}$. We have $\tau(w)=(w_1,w_2, \ldots,
w_k,s_E(w), \ldots )$ and
\[
\sigma_e(\tau(w))=(r_E(e),r_E(e), \ldots, ew_1w_2\cdots
w_m,w_{m+1}, \ldots, w_k, s_E(w), \ldots ),
\]
where $ew_1w_2\cdots w_m$ is in the $m$th position. Note that
$|w_m|=k_m n_{m-1} < n_m-n_{m-1}$ (as in (\ref{wi})), and thus
\begin{eqnarray*}
|ew_1w_2\cdots w_m| &=& 1+(n_1-1)+(n_2-n_1)+ \cdots \\ & & \cdots
+(n_{m-1}-n_{m-2})+k_mn_{m-1} \\ &=& n_{m-1}+k_mn_{m-1} < n_m.
\end{eqnarray*}
This shows that $ew=r_E(e)r_E(e) \cdots (ew_1\cdots w_m) w_{m+1}
\cdots w_k$ is the decomposition of $ew$ as in (\ref{decomp}). It
follows that $\tau(ew)=\sigma_e(\tau(w))$, and this establishes
(i).

For (ii), note first that, since the topology on $X$ is the
product topology and each $X_i$ is a finite set, every subset of
$X$ that is defined by conditions involving only finitely many
coordinates is both closed and open. Thus every subset of $Y$
defined by such conditions is closed and open in the relative
topology of $Y$. This shows that $D_e$ is closed and open in $Y$.

For every $m\in \mathbb{N}\cup \{\infty\}$ we write \[D_m=\{y\in
Y: r_E(y_1)=s_E(e) \;\;\mbox{and}\;\; i(y)=m \}.\] Then each $D_m$
with $m<\infty$ is closed and open in $Y$ and the set
$D_{\infty}=D_e \backslash \cup_{m<\infty} D_m $ is a closed set
in $Y$. We also write $R_m$ ($m\in \mathbb{N}$) for
\[R_m=\big\{y=(y_1,y_2,\ldots )\in Y :
 y_i=r_E(e)
\;\mbox{for all}\;\; i<m \;\]\[\;\mbox{and}\;\; y_m=ew'  \;
\mbox{for some}\; |w'|\equiv -1 (mod \; n_{m-1}) \big\}
\] and
\[R_{\infty}=\{(r_E(e),r_E(e),\ldots ) \}. \]
Then each $R_m$ (with $m<\infty$) is open and closed in $Y$ and
$R_{\infty}$ is closed. Also $R_e$ is the disjoint union of all
the $R_m$'s.

Fix $m<\infty$ and define the restriction $\sigma_m =
\sigma_e|_{D_m}$. It is easy to see that $\sigma_m$ is a
homeomorphism from $D_m$ onto $R_m$; in fact, it is injective and
involves a change in only finitely many coordinates. We also know
that $\sigma_e$ maps $D_{\infty}$ onto the (one-point) set
$R_{\infty}$. Thus $\sigma_e$ maps $D_e$ onto $R_e$ and its
restriction to the complement of $D_{\infty}$ is continuous.

Suppose $\{x^n\}$ is a sequence in $R_e$ converging to some $y\in
Y$. If $y$ is not in $R_e$ then, for every $m<\infty$, only
finitely many elements of the sequence lie in $R_m$. Thus, for
every $m<\infty$ we can find some $K_m$ such that for every
$k>K_m$, $x^k$ is not in $R_i$ for $i\leq m$. Thus, for $k>K_m$,
$x^k_i=r_E(e)$ for all $i\leq m$. It follows that the limit, $y$,
is equal to $(r_E(e),r_E(e), \ldots )$ and, thus, lies in $R_e$.
Therefore $R_e$ is closed in $Y$.

It is left to show that $\sigma_e$ is continuous. In fact, it is
left to consider sequences $\{z^k\}$ in $D_e$ converging to some
$z\in D_{\infty}$. But then the sequence $\{\sigma_e(z^k)\}$ lies
in $R_e$. Since $Y$ is compact (and the topology is metric) we can
find a converging subsequence. As the argument above shows, the
limit will be in $R_e$ and, in fact, it will be in $R_{\infty}$
(since $z\in D_{\infty}$). Since $R_{\infty}$ has only one point
and this point is the image of $z$ under $\sigma_e$, we are done.
 \bx

We now use the notation set above to introduce the topological
graph $E(\infty)$ that will play an important role in studying the
algebra $\B_E(\{n_k\})$. Topological graphs were introduced by the
second author and Muhly in \cite{MS_Morita} (see Example 5.4
there) and studied further by Katsura (\cite{Ka1,Ka2,Ka3}) and by
Muhly and Tomforde (\cite{MT}). They are also referred to as
``continuous graphs" or ``topological quivers". (More precisely,
the definition in \cite{MS_Morita} and in \cite{MT} of a
topological quiver is somewhat more general as it requires only
that the map $s$ is continuous and open and not necessarily a
local homeomorphism).

\begin{defn}
A {\it topological graph} is given by a quadruple $F=(F^0,F^1,
s_F,r_F)$ where $F^0$, $F^1$ are locally compact spaces, $s_F: F^1
\rightarrow F^0$ is a local homeomorphism and $r_F:F^1 \rightarrow
F^0$ is a continuous map.
\end{defn}

To a topological graph $F$ one associates a graph $C^*$-algebra,
written $\O(F)$ (\cite{Ka1} and \cite{MT}). We will not go into
the details of the definition of $\O(F)$ but just note that it
generalizes $\O(F)$ for finite graphs and it is the Cuntz-Pimsner
$C^*$-algebra associated with a $C^*$-correspondence constructed
from the graph. For the graph $E(\infty)$ that we define below,
the $C^*$-correspondence will be discussed later (in
$\S$~\ref{S:Ktheory}).

Now we define the topological graph $E(\infty)$ as follows. Let
\[E(\infty)^0=Y,\]\[ E(\infty)^1=\{(e,\omega)\in E^1 \times Y : \omega\in D_e
\}\] and for all $(e,\omega)\in E(\infty)^1$ put \[
s_{E(\infty)}(e,\omega)=\omega \qand
r_{E(\infty)}(e,\omega)=\sigma_e(\omega).
\] Note that both $E(\infty)^0$ and $E(\infty)^1$ are compact
spaces, the map $s_{E(\infty)}$ is a local homeomorphism (since
its restriction to each $\{e\}\times D_e$ is a homeomorphism onto
$D_e$) and $r_{E(\infty)}$ is continuous (since each $\sigma_e$
is).

Recall now that, given $n,k\in \mathbb{N}$, we defined a regular
factor map $m$ from $E(nk)$ to $E(n)$ (see the discussion that
precedes  Lemma~\ref{m}). With the sequence $\{n_k\}$ as above, we
have a regular factor map from $E(n_k)$ to $E(n_{k-1})$ and we
denote it by $m_{k-1,k}$. We also define, for every $k\in
\mathbb{N}$, a pair $m_k=(m_k^0,m_k^1)$ of maps where $m_k^0$ maps
$E(\infty)^0$ onto $E(n_k)^0$ is defined by
\[ m_k^0(w_1,w_2, \ldots )=w_1w_2\cdots w_k \in E^{<n_k}=E(n_k)^0
\] and $m_k^1$ maps $E(\infty)^1$ onto $E(n_k)^1$ defined by
\[ m_k^1(e,\omega)=(e,m_k^0(\omega))\;,\;\; (e,\omega)\in E(\infty)^1 .\]

These maps are continuous and satisfy
\[ m^0_{k-1,k} \circ m_k^0 = m_{k-1}^0 \] and
\[ m^1_{k-1,k} \circ m_k^1 = m_{k-1}^1 \] for all $k\in
\mathbb{N}$. Also, it is easy to check that, given a sequence
$\{w^{(k)}\}$ where $w^{(k)}\in E(n_k)^0$ for all $k$ and
$m_{k-1,k}(w^{(k)})=w^{(k-1)}$, there is a unique $w\in E(\infty)$
with $m_k(w)=w^{(k)}$ for all $k$. Similar considerations work for
the edge spaces. We also have
$s_{E(n_k)}(m_k^1(e,w))=m_k^0(s_{E(\infty)}(e,w))$ for every
$(e,w)\in E(\infty)^1$ and a similar equality holds for the range
maps. In fact, we see that $E(\infty)$ is the {\it projective
limit} \cite{Ka2} of the projective system defined by the graphs
$E(n_k)$ and the maps $m_{k-1,k}$. (See \cite[$\S$~4]{Ka2}). Also
note that this projective system is surjective (i.e. each
$m_{k-1,k}^0$ is a surjective map).

We may now prove the following.

\begin{thm}\label{limit}
Let $E$ be a finite graph with no sinks and no sources and let
$\{n_k\}$ be an increasing sequence of positive integers with each
$n_k$ dividing $n_{k+1}$. Then the algebra $\B_E(\{n_k\})$  is
$*$-isomorphic to the Cuntz-Pimsner $C^*$-algebra $\O(E(\infty))$.
\end{thm}

\Prf By Theorem~4.13 of \cite{Ka2}, $\O(E(\infty))$ is isomorphic
to the direct limit of the algebras $\O(E(n_k))$ with respect to
the maps $\mu_{m_{k-1,k}}$. Corollary~\ref{jm} shows that the maps
$\pi_{n_k}$ of Theorem~\ref{Bn} can be used to get an isomorphism
of this direct limit and the direct limit of the algebras
$\B(n_k)$ with respect to the maps $j_{n_k,n_{k-1}}$. This
concludes the proof since the latter algebra is $\B_E(\{n_k\})
=\lim_k \B(n_k)$.
 \bx

\begin{rem}\label{limAn}
One can also construct a topological graph $E[\infty]$ satisfying
$\A(n)\cong \O(E[n])$. It is the projective limit of the graphs
$E[n_k]$ (see Theorem~\ref{A(n)}) with respect to the factor maps
defined as in Remark~\ref{factorq}.
\end{rem}

\begin{rem}\label{Etilderemark}
Our no source assumption on $E$ was made to clarify the
presentation. In general though, one can consider the graph
$\tilde{E}$ defined by $E$ as follows:
\[
\tilde{E}^0=\{ v\in E^0 : |\{w\in E^* : r_E(w)=v \}|=\infty \},
\]
\[
\tilde{E}^1=E^1 \cap s_E^{-1}(\tilde{E}^0) \cap
r_E^{-1}(\tilde{E}^0),
\]
\[
s_{\tilde{E}}=s_E|\tilde{E}^1 \qand r_{\tilde{E}}=r_E|\tilde{E}^1.
\]
It is easy to check that $\tilde{E}$ has no sources and no sinks
(provided $E$ has no sinks). Also, if $E$ has no sources and no
sinks then $\tilde{E}=E$.

Then for an arbitrary finite graph $E$ with no sinks (possibly
with sources), $E$ may be replaced by $\td{E}$ in Theorem~\ref{Bn}
and Lemma~\ref{m}. Further, in Theorem~\ref{limit},
$\O(E(\infty))$ may be replaced by $\O(\td{E}(\infty))$. The
following result is an immediate consequence of this generalized
version of the previous theorem.
\end{rem}

\begin{cor}
If $E$ and $F$ are two finite graphs with no sinks and $\tilde{E}$
is isomorphic to $\tilde{F}$ (in particular, if $F=\tilde{E}$)
then $\B_E(\{n_k\})$ and $\B_F(\{n_k\})$ are isomorphic.
\end{cor}

\section{Example}\label{S:example}

Let us denote by $C_j$ the directed graph which is a single simple
loop (or `cycle') with $j$ vertices. In this section we shall
discuss the algebras $\B_{C_j}(\{n_k\})$. Note that the algebra
$\B_{C_1}(\{n_k\})$ is the classical Bunce-Deddens algebra
associated with the sequence $\{n_k\}$.

Fix a positive integer $j$. Write $v_1,v_2, \ldots , v_j$ for the
(distinct) vertices of $C_j$ and $e_1,e_2,\ldots ,e_j$ for its
edges where $s(e_i)=v_i$, $r(e_i)=v_{i+1}$ if $i<j$ and
$r(e_j)=v_1$.

Given a positive integer $n$,  write $p$ for the least common
multiple of $n$ and $j$ and $l$ for their greatest common divisor
(so that $lp=jn$). We write  $n=qj+r$ for the the division of $n$
by $j$, where $q,r$ are integers and $0\leq r <j$. Then
$\gcd(j,r)=\gcd(j,n)=l$ and, considering the equivalence relation
given by addition modulo $r$ on $\{1,2, \ldots ,j\}$, there are
$l$ equivalence classes (each containing $j/l=p/n$ elements). Let
$\Omega$ be a fixed set of representatives, one for each
equivalence class.

The graphs $C_j(n)$ below are the graphs $E(n)$ (of
$\S$~\ref{S:B(n)}) with $E=C_j$.

\begin{lem}\label{Cm}
Let $j,n$ be two positive integers. Then $C_j(n)$ is a disjoint
union of $l$ loops, each of length $p$. In fact,
\begin{equation}\label{decomposition}
 C_j(n)=\sqcup_{i\in \Omega} C_j(n)^{(i)}\end{equation}
where $C_j(n)^{(i)}$ is the loop that contains the vertex $v_i$.

Thus, for each $i\in \Omega$ there is an isomorphism
$\phi_{i,n}=(\phi_{i,n}^0,\phi_{i,n}^1)$
 from $C_1(p)$ to $C_j(n)^{(i)}$  and we can write
 \[ C_j(n)\cong C_1(p) \sqcup C_1(p) \sqcup \ldots \sqcup C_1(p),
 \] where the right hand side is a disjoint union of $l$ simple
 loops.
\end{lem}

\Prf For every vertex $v\in E=C_j$, there is a single edge ending
at $v$ and a single edge emanating from it. It follows from this
that the same holds for $C_j(n)$. Thus $C_j(n)$ is a disjoint
union of loops.

Now fix $i\in \Omega$. Start with the vertex $v_i$ in $C_j(n)^0$,
recalling that $C_j^0\subseteq C_j(n)^0$. Travelling along  the
edges in $C_j(n)^1$ we will, after $n-1$ `moves' reach a vertex
$w$ with $|w|=n-1$ and $s_E(w)=v_i$ (there is only one such $w$).
From there the only way to proceed is
 along the edge in $C_j(n)^1$ whose source is $w$. This edge is
 $(e_{i\oplus (n-1)},w)=(e_{i \oplus (r-1)},w)$, where we write $\oplus$
 for addition modulo  $j$. Its range is
 $$r_{C_j(n)}(e_{i\oplus (r-1)},w)=r_{C_j}(e_{i\oplus (r-1)})=v_{i \oplus
 r}.$$
 Thus, after `moving along'  $n$ edges (starting at $v_i$)  we reach the vertex
 $v_{i \oplus r}$. Travelling along $n$ more edges we reach $v_{2r\oplus i}$ and so
on until we get back to $v_i$. Clearly, $\{i,r\oplus i, 2r \oplus
i, \ldots \}$ is one of the equivalence classes mentioned above.
In fact, it is the equivalence class whose
 representative (in $\Omega$) is $i$ and it
contains $p/n$ elements. Thus, this loop contains $(p/n) n=p$
edges (and vertices) and we denote it by $C_j(n)^{(i)}$.
 Since this argument holds
for every loop,
 it shows that each loop contains $p$ edges,
completing the proof of the first statement of the lemma.

The last statement of the lemma is clear since all simple loops of
length $p$ are isomorphic.
 \bx

\begin{rem}\label{loop}
As was mentioned in the lemma, each loop $C_j(n)^{(i)}$ is
isomorphic to the graph $C_1(p)$. In fact, there are $p$ different
directed graph isomorphisms from $C_1(p)$ to $C_j(n)^{(i)}$.
 We wish to fix
one and we do it as follows. If $C_1^0=\{v\}$ and $C_1^1=\{e\}$,
then a vertex in $C_1(p)^0$ is of the form $v$ or $eee\cdots e$
(with no more than $p-1$ $e$'s). For each $i\in \Omega$ we fix the
only isomorphism from $C_1(p)$ to $C_j(n)^{(i)}$ that sends the
vertex $v$ (in $C_1(p)$)
 to the vertex
$v_i$ (in the loop $C_j(n)^{(i)}$). From now on, when we write
$\phi_{i,n}=(\phi_{i,n}^0,\phi_{i,n}^1)$, we refer to this
isomorphism.
\end{rem}

\begin{lem}\label{mm}
Let $j,n,k$ be positive integers such that $\gcd(j,nk)=\gcd(j,n)$
(and write $l$ for this number and $p$ for $nj/l$). Let
$m=(m^0,m^1)$ be the factor map from $C_j(nk)$ to $C_j(n)$ as in
Lemma~\ref{m}. Then, for every $i\in \Omega$, $m$ maps
$C_j(nk)^{(i)}$ into $C_j(n)^{(i)}$. Moreover, writing
$\phi_{i,nk}$ and $\phi_{i,n}$ for the isomorphisms in
Remark~\ref{loop} associated, respectively, with $C_j(nk)^{(i)}$
and $C_j(n)^{(i)}$, the map $(\phi_{i,n})^{-1} \circ m \circ
\phi_{i,nk}$ is the factor map from $C_1(pk)$ to $C_1(p)$ as in
Lemma~\ref{m}.
\end{lem}

\Prf The factor map $m$, as defined in the discussion that
proceeds Lemma~\ref{m}, maps $E(nk)$ to $E(n)$ and fixes the
vertices in $E^0$. (Recall that $E^0$ is contained in both
$E(n)^0$ and in $E(nk)^0$.) Thus, when $E=C_j$, it fixes the
vertices $v_1 ,\ldots,v_j$ and, in particular, it fixes each $v_i$
for $i \in \Omega$. So fix such an $i$ and write the vertices of
the loop $C_j(nk)^{(i)}$ as $\{u_1=v_i, u_2, u_3, \ldots
,u_{pk}\}$ and the vertices of $C_j(n)^{(i)}$ as
$\{z_1=v_i,z_2,z_3, \ldots, z_p\}$. Also write $f_q $ for the edge
in $C_j(nk)^{(i)}$ emanating from $u_q$ and ending at $u_{q+1}$
and, similarly write $g_q$ for the edge of $C_j(n)^{(i)}$ starting
at $z_q$ and ending at $z_{q+1}$ (with $f_{pk}$
 and $g_{p}$ defined in the obvious way).

As $m^0(u_1)=z_1$, it follows from (i) of Definition~\ref{factor}
that $s(m^1(f_1))=m^0(s(f_1))=m^0(u_1)=z_1$ and, consequently,
$m^1(f_1)=g_1$. Using Definition~\ref{factor} again, we get
$m^0(u_2)=m^0(r(f_1))=r(m^1(f_1))=r(g_1)=z_2$. (Here we used $r,s$
to denote the range and source maps for both graphs but that
should cause no confusion.) Continuing in this way we see that $m$
maps $C_j(nk)^{(i)}$ onto $C_j(n)^{(i)}$. In fact, the image of
$m$ `travels' along the smaller loop $k$ times.

This argument shows, in fact, that there is a unique factor map
from a loop of length $pk$ to a loop of length $p$, provided we
require that a chosen vertex in the first loop will be mapped to a
chosen one in the second. Since $C_1(pk)$ and $C_1(p)$ are such
loops and the map $(\phi_{i,n})^{-1} \circ m \circ \phi_{i,nk}$ is
a factor map from $C_1(pk)$ to $C_1(p)$ that maps the vertex $v$
(in $C_1(pk)$) to $v$ (in $C_1(p)$) , it is the unique factor map
that does it. It follows that it equals the factor map of
Lemma~\ref{m} (with $C_1$ in place of $E$ and $p$ in place of
$n$).
 \bx

\begin{cor}\label{BCj}
Let $j,n,k$ be positive integers such that $\gcd(j,nk)=\gcd(j,n)$
(and write $l$ for this number and $p$ for $nj/l$). Then there are
$*$-isomorphisms
\[  \Phi_n : \B_{C_j}(n) \rightarrow \B_{C_1}(p) \oplus \B_{C_1}(p) \oplus
\cdots \oplus \B_{C_1}(p) \] and \[
 \Phi_{nk} : \B_{C_j}(nk) \rightarrow \B_{C_1}(pk) \oplus \B_{C_1}(pk) \oplus
\cdots \oplus \B_{C_1}(pk) \] such that, for every $i\in \Omega$,
\[j_{nk,n}^{C_j}|_{C_j(n)^{(i)}} =(\Phi_{nk})^{-1} \circ j_{pk,p}^{C_1} \circ
\Phi_n |_{C_j(n)^{(i)}} \] where $j_{nk,n}^{C_j}$ and
$j_{pk,p}^{C_1}$ are the maps defined in (\ref{j}) associated with
the graphs $C_j$ and $C_1$ respectively. Hence
\[j_{nk,n}^{C_j}=(\Phi_{nk})^{-1} \circ \Big(\sum_{i\in \Omega} \oplus j_{pk,p}^{C_1}\Big) \circ
\Phi_n.  \]
\end{cor}

\Prf The isomorphisms $\Phi_n$ and $\Phi_{nk}$ are the ones
implemented by the graph-isomorphisms $\sum_i \phi_{i,n}$ and
$\sum_i \phi_{i,nk}$ respectively. (These maps are defined in
Remark~\ref{loop}.) Since, by Corollary~\ref{jm}, the maps
$j_{nk,n}^{C_j}$ and $j_{nk,n}^{C_1}$ are the ones implemented by
the corresponding factor maps, the result follows from
Lemma~\ref{mm}.
 \bx

\begin{thm}\label{bdclass}
For a positive integer $j$ and a sequence of positive integers
$\{n_k\}$ as above, the $C^*$-algebra $\B_{C_j}(\{n_k\})$ is
$*$-isomorphic to the direct sum of $l$ copies of the classical
Bunce-Deddens algebra $\B_{C_1}(\{p_k\})$
 where $l=\max_k \gcd(j,n_k)$ and $p_k=jn_k/l$.

 It follows that
 \[\B_{C_j}(\{n_k\}) \cong \B_{C_{j'}}(\{n'_k\}) \]
 if and only if $\max_k \gcd(j,n_k)=\max_k \gcd(j',n'_k)$ ($=l$)
 and the supernatural numbers associated with $\{jn_k/l\}$ and
 with $\{j'n'_k/l\}$ coincide.
 \end{thm}

\Prf The sequence $\{\gcd(j,n_k)\}_{k=1}^{\infty}$ is a non
decreasing sequence of positive integers that are smaller or equal
to $j$. Thus, for some $k_0$, $\gcd(j,n_k)=l$ whenever $k\geq
k_0$. Since we are interested in the limit algebra, we can, and
shall, assume that $\gcd(j,n_k)=l$ for all $k$.

Thus, we can use Corollary~\ref{BCj} and the fact that
\[\B_{C_j}(\{n_k\})=\lim_k \Big(\B_{C_j}(n_k),j^{C_j}_{n_{k+1},n_k}\Big) \]
and
\[\B_{C_1}(\{p_k\})=\lim_k
\Big(\B_{C_1}(p_k),j^{C_1}_{p_{k+1},p_k}\Big)\] to conclude that
the family of $*$-isomorphisms $\{\Phi_{n_k} \}$ (defined in
Corollary~\ref{BCj}) defines a $*$-isomorphism from
$\B_{C_j}(\{n_k\})$ onto the direct sum of $l$ copies of
$\B_{C_1}(\{p_k\})$, completing the proof of the first statement
of the theorem.

Now assume that \[\B_{C_j}(\{n_k\}) \cong \B_{C_{j'}}(\{n'_k\}).\]
Note that the $C^*$-algebra $\B_{C_1}(\{p_k\})$ is the classical
Bunce-Deddens algebra associated with the sequence $\{p_k\}$. It
is known to be simple (\cite[Theorem V.3.3]{D}) and, thus, the
center of $\B_{C_j}(\{n_k\}) $ is of dimension
$l=\max_k\gcd(j,n_k)$ and is generated by an orthogonal family of
$l$ projections whose sum is $I$. It then follows that
$\max_k\gcd(j,n_k)=\max_k\gcd(j',n'_k)$ (since the centers of the
two algebras are isomorphic). Also, if $q$ is one of these central
projections in $\B_{C_j}(\{n_k\})$ and it is mapped by the
isomorphism to the central projection $q'$ in the other algebra,
then the isomorphism maps $q\,\B_{C_j}(\{n_k\})\, q$ (which is
isomorphic to $\B_{C_1}(\{jn_k/l\})$) onto the algebra
$q'\,\B_{C_j'}(\{n'_k\})\,q'$ (which is isomorphic to
$\B_{C_1}(\{j'n'_k/l\})$). It follows from \cite[Theorem V.3.5]{D}
that the two supernatural numbers coincide. Theorem V.3.5 of
\cite{D}, together with the first statement of the theorem, proves
the other direction.
 \bx

\begin{cor}\label{simpleloop}
The algebra $\B_{C_j}(\{n_k\})$ is simple if and only if for every
$k\geq 1$, we have $\gcd(j,n_k)=1$ .
\end{cor}

\begin{rem}
We expect that the classification result of \cite{Kinductive},
which generalizes the Bunce-Deddens supernatural number
classification to the Cuntz case, could be used to extend
Theorem~\ref{bdclass} to a broader class of graphs $E$. For the
sake of succinctness, however, we shall not pursue this further
here. More generally, we wonder for what graphs $E$ could a
classification theorem along the lines of Theorem~\ref{bdclass} be
proved?
\end{rem}

\section{K-Theory}\label{S:Ktheory}

In this section we derive the $K$-groups of the algebra
$\B_E(\{n_k\})$, where again $E$ is a finite graph with no sinks
and no sources and $\{n_k\}$ is a sequence as above. We start with
the $K$-theory of $C(Y)$.

\begin{lem}\label{KofA}
Let $Y$ be the $\{n_k\}$-compactification of $E^*$. Then
\[ K_0(C(Y))\cong C(Y,\mathbb{Z}) \]
and
\[ K_1(C(Y))=\{0\} .\]
\end{lem}

\Prf For every $k \in \mathbb{N}$ write $\C_k$ for the subalgebra
of all functions $f$ in $C(Y)$ with the property that $f(y)=f(z)$
whenever $y_i=z_i$ for all $i\leq k$. There is a $*$-isomorphism
$\rho_k:C(E(n_k)^0)\rightarrow \C_k$ given by
$\rho_k(g)(y)=g(y_1y_2 \cdots y_k)$. If $\iota_{k+1,k}$ is the
inclusion map of $\C_k$ into $\C_{k+1}$ then the map
$\rho_{k+1}^{-1}\circ \iota_k \circ \rho_k$ is equal to the map
$\mu_{n_{k+1},n_k}^0$ defined above. Note that $\cup_k \C_k$ is a
dense subalgebra of $C(Y)$ (by the Stone-Weierstrass Theorem).
Thus
\[ C(Y)=\lim_{k} \big(C(E(n_k)^0), (\mu_{n_{k+1},n_k})^0\big) .\]

Fix $f\in C(Y)$ with values in $\mathbb{Z}$. For $0<\epsilon <1/2$
we can, by the above, find $k$ and $g\in \C_k$ with $\Vert f-g
\Vert <\epsilon$. Let $\psi:\cup_{n\in\mathbb{Z}}
(n-\epsilon,n+\epsilon) \rightarrow \mathbb{Z}$ be defined by
$\psi |_{(n-\epsilon,n+\epsilon)} \equiv n$. Then $\psi$ is
continuous and so is the function $g':=\psi \circ g$. But clearly
$g'\in \C_k$ and $f=g'$. Thus $f\in \C_k$.

 This shows that  $$\big\{f\in
C(Y,\mathbb{Z}):\; {\rm for\, some} \, k,\, f\in \C_k\big\}=
C(Y,\mathbb{Z}).$$ Using the notation $((\mu_{n_{k+1},n_k})^0)_*$
for the restriction of this map to $\mathbb{Z}$-valued functions
in $C(E(n_k)^0)$, we get
\[ C(Y,\mathbb{Z})=\lim_{k} \big(C(E(n_k)^0,\mathbb{Z}), ((\mu_{n_{k+1},n_k})^0)_* \big).\]
Since $K_0(C(E(n_k)^0))$ is isomorphic to $C(E(n_k)^0,\mathbb{Z})$
(recall that $E(n_k)^0$ is a finite set) and
$((\mu_{n_{k+1},n_k})^0)_*$ is the map induced from
$(\mu_{n_{k+1},n_k})^0$ on the $K_0$ groups, we find that
\begin{equation}\label{K0cy}
K_0(C(Y))\cong \lim_k \big(C(E(n_k)^0,\mathbb{Z}),
((\mu_{n_{k+1},n_k})^0)_* \big)\cong
C(Y,\mathbb{Z}).\end{equation} Since $K_1(C(E(n_k)^0))=\{0\}$ for
each $k$, the second statement of the lemma also follows.
 \bx

Given the topological graph $E(\infty)$, one can associate with it
a $C^*$-correspondence $Z$ over $A=C(Y)$ as follows. (See also
\cite{MS1},\cite{MT} and \cite{Ka1}). On the space
$C(E(\infty)^1)$ one can define a (right) $C(Y)$-module structure
by setting
\[ (\psi f)(e,w)=\psi(e,w) f(w),\;\; \psi \in C(E(\infty)^1),\;f
\in C(Y)\] and a $C(Y)$-valued inner product by
\[ \langle \psi_1,\psi_2 \rangle (y)=\sum_{(e,y)\in E(\infty)^1}
\overline{\psi_1(y)} \psi_2(y) .\] This makes  $C(E(\infty)^1)$
into a Hilbert $C^*$-module over $A=C(Y)$. (Note that the sum
defining the inner product is finite so that we don't have to
replace $C(E(\infty)^1)$ by $C_d(E(\infty)^1)$ as in \cite{Ka1} or
by a completion of $C_c(E(\infty)^1)$ as in \cite{MT}. Note also
that, using \cite[Lemma 1.8]{Ka1}, this module is already a Banach
space with respect to the norm defined by the inner product.)

To make this module into a correspondence one defines, for $f\in
C(Y), \psi \in C(E(\infty)^1)$,
$(f\psi)(e,w)=f(\sigma_e(w))\psi(e,w).$ This defines the
correspondence associated with this graph. It will be convenient,
however, to write it in a slightly different way. For this, note
first that, for $e\in E^1$, $C(D_e)$ is a Hilbert $C^*$-module
over $A=C(Y)$ and can be made into a $C^*$-correspondence by
defining the left action using $\sigma_e$
\[ (f\cdot g)(y)=f(\sigma_e(y))g(y),\;\;f\in C(Y),\, g\in C(D_e).\]

Now we let $Z$ be the correspondence
\[ Z=\oplus_{e\in E^1} C(D_e).\]
We write $\phi_Z$ for the left action; i.e.,
\[ \phi_Z(f) (\oplus g_e) = \oplus (f\circ \sigma_e) g_e .\]

Given $\psi \in C(E(\infty)^1)$ and $e\in E^1$, write  $\psi_e \in
C(D_e)$ for the function $\psi_e(w)=\psi(e,w)$. Then it is easy to
check that the map
\[ \psi \mapsto \oplus \psi_e \]
is an isomorphism of correspondences from $C(E(\infty)^1)$ onto
$Z$. Thus, we can write $Z$ for the correspondence associated with
the graph $E(\infty)$.

In order to state the next result, one should note that $Z$ is a
finitely generated Hilbert $C^*$-module over $A$ and the triple
$(Z,\phi_Z,0)$ defines an element in $KK(A,A)$. As such, it
defines a map on $K_0(A)$ (into itself), written $[Z]$. In fact, a
general element of $K_0(A)$ can be written as a difference
$[\mathcal{E}_1]-[\mathcal{E}_2]$ for finitely generated
projective modules
 $\mathcal{E}_i$ over $A$ and the map $[Z]$, defined by
 $(Z,\phi_Z,0)$,
 will map it into the element $[\mathcal{E}_1\otimes_A Z]-[\mathcal{E}_2
 \otimes_A Z]$.

 Using Lemma~\ref{KofA}, it follows that $[Z]$ induces a map on
 $C(Y,\mathbb{Z})$. To see how  this map is defined we first need
 the following discussion.

Given a (finite) subset $B\subseteq E(n_k)^0$ for some $k\geq 1$,
and, given some $j\geq k$, we form
\[ B(j)=\{w\in E(n_j)^0 : w(n_k)\in B \} \] and
\[ B(\infty)=\big\{y=(y_1,y_2, \ldots )\in Y : y_1y_2 \cdots y_k \in B
\big\}.\] Then $B(\infty)$ is a subset of $Y$ that is both closed
and open. In fact, every subset of $Y$ that is closed and open is
$B(\infty)$ for some $k$ and some subset $B$ of $E(n_k)^0$.

For such $B$ write $\chi_{B(\infty)}$ for the characteristic
function of $B(\infty)$. Then $\chi_{B(\infty)}\in
C(Y,\mathbb{Z})$. Write
\[J_B=\big\{g \in C(Y) : g(y)=0 ,\;y\in Y\backslash B(\infty) \big\}.\]
Then $J_B$ is a finitely generated projective $C(Y)$-module. Thus,
it defines an element $[J_B]$ in $K_0(C(Y))$. The function in
$C(Y,\mathbb{Z})$ associated with this element via the isomorphism
of Lemma~\ref{KofA} is $\chi_{B(\infty)}$. To see this, one can
write $J_B$ as a direct limit of $J_B \cap \C_j$ and note that
$J_B \cap \C_j$ defines the element in $K_0(\C_j)\cong
C(E(n_j)^0,\mathbb{Z})$ which is the characteristic function of
$B(j)$. For $j\geq k$ the image of this function, under the
embedding of $C(E(n_j)^0,\mathbb{Z})$ into $C(Y,\mathbb{Z})$ given
by the direct limit (\ref{K0cy}), is the characteristic function
of $B(\infty)$.

For $B$ as above and $e\in E^0$, consider the set
$\sigma_e^{-1}(B(\infty))$. It is also a closed and open subset of
$Y$ and, thus, is equal to $C(\infty)$ for some $k$ and
$C\subseteq E(n_k)^0$. We have $J_C=\{g \in C(Y) : g(y)=0 ,\;y\in
 Y\backslash C(\infty) \}=\{g \in C(Y) : g(y)=0 ,\;y \in
 Y\backslash \sigma_e^{-1}(B(\infty))\}=\{f\circ \sigma_e : f\in
 J_B \}$. We write $J_B \circ \sigma_e$ for this space. (Here, and
 below, the function $f\circ \sigma_e$ is assumed to vanish
 outside $D_e$).

For $B$ as above we now consider $J_B \otimes_A Z$. It is easy to
see that this Hilbert $C^*$-module is isomorphic to
$\phi_Z(J_B)Z=\oplus_e (J_B \circ \sigma_e)$. It follows that
$[Z]([J_B])=\sum [J_B \circ \sigma_e]$ and, viewing
 $[Z]$ as a map of $C(Y,\mathbb{Z})$ (via the isomorphism
of Lemma~\ref{KofA}), we get
\[ [Z](\chi_{B(\infty)})=\sum \chi_{B(\infty)} \circ \sigma_e .\]
Since every closed and open set in $Y$ is of the form $B(\infty)$,
these characteristic functions span $C(Y,\mathbb{Z})$. Thus
\begin{equation}\label{Z}
[Z](f)=\sum_e (f \circ \sigma_e)
\end{equation}
where $f\in C(Y,\mathbb{Z})$ and $(f \circ \sigma_e)(y)$ is
understood to be $0$ if $y$ is not in $D_e$.

 Applying  a result of Katsura (\cite[Corollary
 6.10]{Ka1}) we get the following. (Note that, in the notation of \cite{Ka1},
 $E(\infty)^0_{rg}=E(\infty)^0$ since
 $r_{E(\infty)}(E(\infty)^1)=E(\infty)^0$ and $E(\infty)^0$ is
 compact.)

\begin{thm}\label{Kgroup} (\cite{Ka1})
Let $Z$ be the correspondence defined above and $[Z]$ be the map
it induces in $K$-theory. Let $t^0$ and $t^1$ be the imbeddings of
$C(E(\infty)^0)$ and of $Z$ into $\mathcal{O}(E(\infty))$
respectively. Then we have the following exact sequence of
$K$-groups: $$\begin{CD} K_0(C_0(E(\infty)^0)) @>>\iota_*-[Z]>
K_0(C_0(E(\infty)^0)) @>>t^0_*>  K_0(\mathcal{O}(E(\infty))) \\
@AAA @. @VVV \\ K_1(\mathcal{O}(E(\infty))) @<t^0_*<<
K_1(C_0(E(\infty)^0)) @<\iota_*-[Z]<< K_1(C_0(E(\infty)^0)).
\end{CD}$$
\end{thm}

\vspace{4mm}

For $f\in C(Y,\mathbb{Z})$ write
\[ \Delta (f)=f - \sum_{e\in E^1} f \circ \sigma_e \]
where $(f\circ \sigma_e)(y)$ is understood to be $0$ if $y$ is not
in $D_e$.

The following theorem is now a direct consequence of
Theorem~\ref{Kgroup}, Lemma~\ref{KofA} and equation (\ref{Z}).

\begin{thm}\label{Kthy}
\[ K_0(\B_E(\{n_k\}))=C(Y,\mathbb{Z})/Im (\Delta) \]
and
\[ K_1(\B_E(\{n_k\}))=Ker (\Delta) .\]
\end{thm}

\begin{rem}
The classification theorem from \cite{Kinductive} shows that the
$K$-groups can be used to distinguish the algebras $\B_E(\{n_k\})$
in the case of single vertex graphs $E$. We expect that
Theorem~\ref{Kthy} could be used to classify the algebras
associated with broader classes of graphs.
\end{rem}

\section{Simplicity}\label{simp}

Simplicity of $C^*$-algebras associated with topological graphs
was characterized in \cite[Theorem 10.2]{MT} and in \cite[Theorem
8.12]{Ka3}. We apply these results to the graph $E(\infty)$. We
first need the following.

\begin{lem}\label{noloop}
The graph $E(\infty)$ contains no loops.
\end{lem}

\Prf Suppose $f_1f_2f_3 \cdots f_k$ is a loop in $F=E(\infty)$ and
let $u^i=s_{F}(f_i)$ for $1\leq i \leq k$. Recall that $u^i\in Y$
and $ u^i_m$ is its $m$th coordinate. We distinguish two cases.

First suppose none of the $u^i$'s lies in $\tau(E^0)$. Then there
is some $N$ such that for $m> N$ and $1\leq i,j \leq k$,
\[ u^i_m=u^j_m . \]

Let \[ g(i)=\sum_{m=1}^N |u^i_1u^i_2 \cdots u^i_N|,\qfor 1\leq i
\leq k\] where $|\cdot|$ is the length of an element of $E^*$.

 Then, for $1\leq i <k$, $g(i+1)=g(i)+1>g(i)$ since
$u^{i+1}=\sigma_e(u^i)$ for some $e\in E^1$. A similar argument
shows that $g(1)>g(k)$, and hence this gives a contradiction.

In the second case suppose one of the $u^i$'s lies in $\tau(E^0)$.
Say, $u^1=\tau(v)$ (for some $v\in E^0$). Write
$f_i=(e_i,\omega_i)\in E(\infty)^1$ (so $e_i\in E^1$ and $\omega_i
\in Y$) and then $u^1=\sigma_{e_{k}}\circ \sigma_{e_{k-1}} \circ
\cdots \circ \sigma_{e_1}(u^1)=\tau(e_{k}e_{k-1} \cdots e_1)$ ,
contradicting the fact that $u^1\in \tau(E^0)$ and $\tau$ is
injective.

Since, in either case, we arrive at a contradiction, $E(\infty)$
contains no loops.
 \bx

Using the notation of \cite{Ka3}, it now follows immediately that
$E(\infty)$ is what Katsura calls a `topologically free graph'
and, in the notation of \cite{MT}, the graph satisfies Condition
(L).

In order to discuss simplicity we need also the notion of
minimality. This is defined in both \cite{Ka3} and \cite{MT}. For
the graph $E(\infty)$ both definitions are easily seen to be
equivalent to the following.

\begin{defn}\label{min}
A subset $B\subseteq E(\infty)^0=Y$ is said to be {\it invariant}
if $\sigma_e(y)\in B$ whenever $y\in B\cap D_e$ and there is some
$f \in E^1$ and $z\in D_f \cap B$ such that $\sigma_f(z)=y$. The
graph $E(\infty)$ is said to be {\it minimal} if there is no
proper, nonempty, closed invariant subset of $Y$.
\end{defn}

The following is a consequence of Lemma~\ref{noloop} and
\cite[Theorem 8.12]{Ka3} or \cite[Theorem 10.2]{MT}.

\begin{thm}\label{simple}
The algebra $\B_E(\{n_k\})$ is simple if and only if $E(\infty)$
is minimal.
\end{thm}

We would like, of course, to have a condition on $E$ and the
sequence $\{n_k\}$ that will be necessary and sufficient for the
minimality of $E(\infty)$. So far, we don't have such a condition
for arbitrary graphs but we present a sufficient condition below.
We shall need the following lemma.

\begin{lem}\label{invar}
Every nonempty closed (with respect to $E(\infty)$) invariant
subset $Y_0 \subseteq Y$ contains an element of the form $\tau(u)$
for some $u\in E^0 $. Moreover, for such $u$, $Y_0$ contains every
$\tau(w)$ for $w\in E^*$ with $s_E(w)=u$.
\end{lem}

\Prf Let $Y_0$ be a closed invariant subset of $Y$. Fix $y\in Y_0$
and write it $y=(y_1,y_2, \ldots )$. We can write $y_1=e_1e_2
\cdots e_j$ for some $e_1,e_2, \ldots ,e_j \in E^1$. Then
$y=\sigma_{e_1}(e_2e_3 \cdots e_j,y_2, \ldots )$ and the element
$z=(e_2e_3 \cdots e_j,y_2, \ldots )$ is the unique one satisfying
$y=\sigma_e(z)$ for some $e\in E^1$. It follows from the
invariance of $Y_0$ that $z\in Y_0$. Continuing in this way and
noting that $y=\sigma_{e_1}\circ \sigma_{e_2} \circ \cdots \circ
\sigma_{e_j}(s_E(e_j)=r_E(y_2),y_2, \ldots )$, we find that
$$y[1]:=(r_E(y_2),y_2, \ldots ) \in Y_0.$$ If $y[1]\in \tau(E^0)$,
we are done. Otherwise, we write $y_2=e'_1e'_2 \cdots e'_l$ for
$e'_1, \ldots e'_l \in E^1$. Note that $y[1]=\sigma_{e'_1} \circ
\cdots \circ\sigma_{e'_l}(r_E(y_3),r_E(y_3),y_3, \ldots )$, and we
conclude from the invariance of $Y_0$ that
$$y[2]:=(r_E(y_3),r_E(y_3),y_3, \ldots ) \in Y_0.$$ Continuing
in this way, we get a sequence $y[k]$ in $Y_0$ with
$$y[k]=(r_E(y_{k+1}), \ldots ,r_E(y_{k+1}),y_{k+1}, y_{k+2},
\ldots ).$$ Since $E$ is a finite graph, one of the vertices, say
$u\in E^0$, will appear infinitely many times in the sequence
$\{r_E(y_{k+1})\}$. So, for some increasing sequence of positive
integers $\{k_m\}$, $r_E(y_{k_m+1})=u$ for every $m$. It follows
that the sequence $ y[k_m]$ converges in $Y$ to $(u,u, \ldots
)=\tau(u).$ Since $Y_0$ is closed, $\tau(u)$ belongs to  $Y_0$.

For the last statement of the lemma, fix $w\in E^*$ with
$s_E(w)=u$ and write it $w=e_1e_2 \cdots e_k$ (with $e_i\in E^1$
and $s_E(e_k)=u$). Then $\tau(w)=\sigma_{e_1}\circ \sigma_{e_2}
\cdots \circ \sigma_{e_k}(\tau(u)) $ and it follows from the
invariance of $Y_0$ that $\tau(w)\in Y_0$.
 \bx

\begin{prop}
If, for every $v,u$ in $E^0$ and every $k\geq 1$, there is some
$w\in E^*$ with $s_E(w)=v$, $r_E(w)=u$ and $|w|$ is a multiple of
$n_k$, then the algebra $\B_E(\{n_k\})$ is simple.
\end{prop}

\Prf Suppose the condition in the hypothesis holds and fix a
closed invariant nonempty subset $Y_0$ of $Y$. We shall show that
$Y_0=Y$. Since $\tau(E^*)$ is dense in $Y$, it suffices to show
that $\tau(E^*)\subseteq Y_0$. From Lemma~\ref{invar} we conclude
that there is some $u\in E^0$ such that $\tau(w)\in Y_0$ whenever
$w\in E^*$ with $s_E(w)=u$.

Now we fix $v\in E^0$ and a positive integer $k$. By assumption it
follows that there is some $w[k]\in E^*$ with $s_E(w[k])=u$,
$r_E(w[k])=v$ and $|w[k]|$ is a multiple of $n_k$. Then $y[k]
:=\tau(w[k])$ has the form $$ y[k]=(v,v, \ldots,
v,w_{k+1},w_{k+2}, \ldots w_{m}, u,u, \ldots )$$ where
$w[k]=vv\cdots vw_{k+1}w_{k+2} \cdots w_{m}$ is the decomposition
of $w[k]$ as in (\ref{decomp}). It follows that $y[k] \rightarrow
(v,v, \ldots)=\tau(v)$, proving that $\tau(v)\in Y_0$. The
argument of the last paragraph of the proof of Lemma~\ref{invar}
shows now that $Y_0$ contains every $\tau(w')$ with $s_E(w')=v$.
Since $v$ is arbitrary, $\tau(E^*)\subseteq Y_0$, and this
completes the proof.
 \bx

Recall that in Corollary~\ref{simpleloop} we showed that
simplicity for the algebras $\B_{C_j}(\{n_k\})$ depends only on
$j$ and $\{n_k\}$. This dependence vanishes when $C_j$ is slightly
adjusted.

\begin{cor}
Let $E$ be a graph that consists of a simple loop with at least
one loop edge attached at some vertex. Then $\B_E(\{n_k\})$ is
simple for every choice of sequence $\{n_k\}$.
\end{cor}

\vspace{0.1in}

{\noindent}{\it Acknowledgements.} We would like to thank Allan
Donsig and Michael Lamoureux for organizing a workshop at the
Banff International Research Station (December 2003) and Gordon
Blower and Stephen Power for organizing a conference in Ambleside,
U.K. (July 2004) that helped foster this collaboration. The first
author was partially supported by an NSERC grant. The second
author was supported by the Fund for the Promotion of Research at
the Technion.




\end{document}